\documentclass[12pt]{amsart}

\def\R{{\mathbb R}}
\def\CA{{\mathcal A}}

\def\acts{\triangleright}

\def\ul#1{\underline{#1}}
\def\lm#1{\lambda^{#1}}
\theoremstyle{remark}
\newtheorem{remark}[]{Remark}

\def\ts{\otimes}

\def\acts{\triangleright}

\title{Rieffel's deformation quantization and isospectral deformations}
\author{Andrzej Sitarz}
\email{Andrzej.Sitarz@th.u-psud.fr}
\thanks{${}^\dagger$\ Supported by Marie Curie Fellowship.}
\address{Laboratoire de Physique Theorique, Universit\'e Paris-Sud, 
Bat. 210, 91405  ORSAY Cedex, France}
\begin{document}
\begin{abstract}
We demonstrate the relation between the isospectral deformation and Rieffel's
deformation  quantization by the action of $\R^d$.
\end{abstract}
\maketitle
\section{Introduction}
The isospectral deformation has been recently introduced by Connes 
and Landi \cite{CoLa}, with the examples of the noncommutative 
3 and 4 spheres, providing new examples of noncommutative spectral
triples \cite{Co2}. However, the construction of similar examples, though
in different context,  can be found in earlier works by Kulish and Mudrov 
\cite{Kul} as well as, more recently  in Paschke Ph.D. thesis \cite{Pas}.

The properties of isospectral deformations has been recently a subject
of much interest \cite{CMDV}. In particular, the symmetries
of isospectral geometries (seen as Hopf algebras acting on the deformed
algebras) appear to be a twist by a Cartan subalgebra of the universal
enveloping algebra of Lie algebras, which were symmetries of the undeformed 
algebras\cite{Sit}, and the constructed spectral triples are symmetric in 
the sense of \cite{PaSi}.

The analysis of symmetries and the generalized construction of the 
deformation, as well as the similarity of the construction with the well-known
noncommutative torus  has led to the connections between the 
{\em isospectral deformations} and Rieffel's \cite{Rie1} deformation quantization 
by the action of $\R^n$. In fact, both the general construction as well as the
corresponding symmetries were already defined in a general setup by Rieffel
(see also \cite{Rie2,Rie3,Rie4}).

In this Letter we show that the isospectral deformations as defined in \cite{CoLa} 
are the special case of the Rieffel's construction \cite{Rie1}. 

\section{The isospectral deformations.}

In the original article \cite{CoLa} one starts with the algebra $\CA$ 
of smooth functions on a manifold, with an isometry group of rank 
$r \geq 2$.  This is equivalent to the statement that the torus 
$T^2$ is the isometry subgroup of the algebra. Now, taking the
elements $t$, which are of $C^\infty$ class with respect to the action 
of the torus group: $ t \mapsto \alpha_s(t)$,  one can obtain their 
decomposition as a norm convergent sum of homogeneous elements 
of a given bidegree. Then, for any homogeneous element one may 
define a deformed product (by a left or right twist). By linearity this 
extends to the linear combinations of homogeneous elements 
which are dense in the algebra we have started with. 

To make correspondence with the Rieffel's deformation 
quantization we shall use the dual picture of isometries as 
used in \cite{Sit}.

\begin{remark}
Let $\CA$ be an algebra and the torus $T^2$ be a subgroup 
of automorphisms of $\CA$. Then for every operator $T$, which 
is of class $C^\infty$ relative to the isometry $\alpha_s$, 
the additive group $\R^2$ acts as follow:

\begin{equation}
 [x_1,x_2] \acts t = \alpha_{e^{2 \pi i x_1}, e^{2 \pi i x_2}}(t), \label{ract} 
\end{equation}

Using the generators of torus symmetries $p_1,p_2$, for instance:
$$ p_1 \acts T = \frac{1}{2\pi i} \left( \frac{d}{dx_1} ( [x_1, x_2] \acts t )
 \right)_{|x=0}, $$
the relation (\ref{ract}) could be rewritten as:
$$ [x_1,x_2] \acts t = e^{2\pi i (x_1 p_1 + x_2 p_2)} \acts t. $$
\end{remark}

\begin{remark}
With respect to the action of $p_1,p_2$ the operators, which
are homogeneous of degree $(n_1,n_2)$, behave like: 
\begin{equation}
\begin{array}{ll}
p_1 \acts t = n_1 t, & p_2 \acts t = n_2 t.
\end{array}
\end{equation}
\end{remark}

\begin{remark}
The product in the algebra $\CA$ can be deformed, first on
elements of given degree, and then extending the deformation
by linearity:
\begin{equation}
a * b = ab \lambda^{n_1^a n_2^b}, \label{deform}
\end{equation}
where $\lambda$ is a complex number such that $|\lambda|=1$.  This 
gives the right twist of \cite{CoLa}.  For future reference we shall denote 
the deformed algebra by $\CA_\lambda$. 
\end{remark}

This deformation, can be extended, in the case 
of a differential manifold and a Lie algebra acting 
on the differential functions,  to all $C^\infty$ functions.
\cite{CoLa}.

It would be useful to introduce a {\em quantization map}:
$$ \CA \ni a  \mapsto \ul{a} \in \CA_\lambda, $$
so that the Eq.\ref{deform} could be rewritten as:
\begin{equation}
\ul{a} * \ul{b} = \ul{ab} \lm{n_1^a n_2^b},  
\label{deform2}
\end{equation}

It will be useful to rewrite the above quantization form in 
a more abstract form, using the action of the generators 
$p_i$. \footnote{We use here a particular form of a more 
general formula derived in \cite{Sit}}. 

\begin{remark}
The quantization map (\ref{deform}) could be written as: $\CA_\Psi$:
\begin{equation}
\ul{m} \left( \ul{a} \ts \ul{b} \right) = 
\ul{m \left( \Psi^{-1} \acts (a \ts b) \right)}, \label{prod}
\end{equation}
where $m$ is the multiplication map 

$$m:\CA \ts \CA \ni a \ts b \to ab \in \CA. $$
and $\ul{m}$ is, similarly,  the multiplication in deformed algebra, 
whereas $\Psi$ is:

\begin{equation}
\Psi_c = \lm{-p_1 \ts p_2},
\end{equation}
\end{remark}

\section{Rieffel's deformation by the action of $\R^n$}

\subsection{General construction}
Suppose we have an algebra $A$ and the action of $V=R^n$ on this 
algebra. and a linear map, $J$, from $V'$, the dual of $V$ to $V$, 
such that it is skew-symmetric $J^T = -J$.  Again, one takes a subalgebra 
of elements, which are $C^\infty$ vectors in $A$ for the action
of $\R^n$. 

To make a direct  correspondence with the above case of 
isospectral deformations we restrict ourselves to $n=2$. In analogy 
with Poisson brackets on the function on a manifold we might define 
a Poisson bracket:

\begin{equation}
P(a,b) = \sum_i  \alpha_{J r_i}(a) \alpha_{p_i}(b), \label{pois}
\end{equation}

where $p_i$ is the basis of $\R^2$ and $r_i$ of its dual. Clearly this 
is independent of the choice of the basis (see \cite{Rie1,Rie2}) and 
makes the algebra $A^\infty$ a strict Poisson *-algebra.

Then using the oscillatory integrals one can define a deformed 
product:

\begin{equation}
a \times_J b = \int_V \int_{V'} \alpha_{Jy}(a) \alpha_x(b) e^{2\pi i (y \cdot x)}, \label{rdef}
\end{equation}
which could be recognized as a deformation quantization 
in the direction of the Poisson structure as defined in (\ref{pois}).

\subsection{Equivalence with isospectral deformations.}

To obtain the deformation as defined in (\ref{deform}) we have 
to modify the expression (\ref{rdef}) by allowing arbitrary (not 
necessarily antisymmetric) operator $J$.

Let us calculate explicitly for homogeneous elements $a$ and 
$b$. We parameterize $V$ (coordinates $x$, basis $p_i$) and $V'$ (coordinates $y$, basis $e_i$),
and  with a particular choice of the map $J$.
$$ J e_1 = 0, \;\;\;\;\; J e_2 = \theta p_1. $$

Then:
$$ 
a \times_J b = \int_V \int_{V'} d^2x d^2y 
e^{2\pi i ( \theta y_2^{} n_1^a )} a e^{2\pi i ( x_1^{} n_1^b +  x_2^{} n_2^b)} b 
e^{ 2\pi i (y_1^{} x_1^{} + y_2^{} x_2^{})} = \ldots 
$$

Calculating further using the standard properties 
of oscillatory integrals we obtain:

$$ 
\begin{aligned}
\ldots = a b \int_V d^2x  \int_{V'} \; d^2y e^{2\pi i ( \theta  n_1^a y_2^{} + n_2^b x_2^{} + y_2^{} x_2^{}) } = \\
\phantom{xxx} = ab e^{2\pi i \theta n_1^a n_2^b },
\end{aligned}
$$

which agrees with the definition for the right twist as defined 
in (\ref{deform}) with $ \lambda = e^{2\pi i \theta}$.

Having shown the relation one might use the results valid 
for the deformation quantization to this particular case, for 
instance to the noncommutative four-spheres. On the other 
hand, using similar arguments as in \cite{CoLa} one will 
obtain a vast family of noncommutative spectral triples and
their symmetries, by constructing the deformation quantization 
as described in \cite{Rie1} of, for instance, spin manifolds.

{\bf Acknowledgements:}
While finishing this paper we have learned of a similar work 
by Joseph Varilly \cite{Var}, more general picture shall be also 
presented in a forthcoming work \cite{CMDV}. 

The author would also like to thank Michel Dubois-Violette, 
J.M.Gracia-Bond\'{\i}a, Gianni Landi,  John Madore, Mario Pachke 
and Harold Steinacker for helpful discussions on many topics.


\begin{thebibliography}{C}
\bibitem{Co} A.Connes, {\em Noncommutative geometry}, Academic Press 1994,
\bibitem{Co2} A.Connes, {\em Noncommutative geometry Year 2000}, arXiv:math.QA/0011193,
\bibitem{CoLa} A.Connes, G.Landi {\em Noncommutative manifolds, the instanton
algebra and isospectral deformations}, arXiv:math.QA/0011194,
\bibitem{CMDV} A.Connes, M. Dubois-Violette, in preparation, 
\bibitem{Kul} P.P.Kulish, A.I.Mudrov, {\em Twist-like geometries on a quantum 
Minkowski space.}, On the occasion of the 65-th birthday of Academician Lyudvig 
Dmitrievich Faddeev, Tr. Mat. Inst. Steklova 226 (1999), Mat. Fiz. Probl. 
Kvantovoi Teor. Polya, 97--111
\bibitem{Pas} M.Paschke, {\em \"Uber nichtkommutative Geometrien, ihre Symmetrien 
und etwas Hochenergiephysik}, Ph.D. thesis, Mainz 2001, to appear, 
\bibitem{PaSi} M.Paschke, A.Sitarz, {\em The geometry of 
noncommutative symmetries},  Acta Phys.Pol. B31, No 11,  (2000)
\bibitem{Rie1} M.Rieffel, {\em Deformation quantization for actions of $R^d$}, 
Memoirs AMS, vol 506, AMS, Providence, (1993),
\bibitem{Rie2} M.Rieffel, {\em Quantization and $C^\ast$ aklgebras}, in: 
$C^\ast$ algebras,1943-1993, 
Contemporary Mathematics 167,  Editor:R.E.Dolan, AMS, Providence, (1994)
\bibitem{Rie3} M.Rieffel, {\em Non-compact quantum groups associated with 
Abelian subgroups}, Comm.Math.Phys. 171, p 181, (1995)
\bibitem{Rie4} M.Rieffel, {\em Deformation quantization for actions of $R^d$}, 
Memoirs AMS, vol 506, AMS, Providence, (1993),
\bibitem{Sit} A.Sitarz, {\em Twists and spectral triples for isospectral 
deformations}, to appear,
\bibitem{Var} J.C.V\'arilly, {\em Quantum symmetry groups of noncommutative spheres}, 
arXiv:math.QA/0102065, 
\end{thebibliography}
\end{document}